\documentclass[11pt]{amsart}
\usepackage{amssymb}
\usepackage{epsfig}
\usepackage{graphicx}
\newtheorem{theorem}{Theorem}
\newtheorem{lemma}[theorem]{Lemma}

\begin{document}

\title{Biquandle longitude invariant of long virtual knots}

\author{Maciej Niebrzydowski}
\address{University of Louisiana at Lafayette,
Lafayette, LA 70504}
\email{mniebrz@gmail.com}

\date{August 28th, 2007}

\thispagestyle{empty}

\begin{abstract}

It is known that the number of biquandle colorings of a long virtual knot diagram, with a fixed color of the initial arc, is a knot invariant.
In this paper we construct a more subtle invariant: a family of biquandle endomorphisms obtained from the set of colorings and longitudinal information. 
 
\end{abstract}

\maketitle

\section{Introduction}

In this paper we will show that keeping track of the order of biquandle elements that appear in the coloring of the long virtual knot diagram allows one to define an invariant that is, at least in some cases, stronger than the total number of colorings, or the number of colorings with a fixed color of the initial arc. Let us begin with some preliminary definitions. 

Virtual knot theory was introduced by Louis Kauffman in \cite{K1}.

\noindent\textbf{Definition.} A \textit{virtual knot} is defined as an equivalence class of 
4-valent plane diagrams, with an extra crossing information, where a new type of crossing (called virtual crossing and denoted by a small circle around a double point) is allowed.

Virtual knot theory is a generalization of classical knot theory; if two classical knots are equivalent under generalized Reidemeister moves, then they are equivalent under classical ones
(see \cite{GPV} for the proof).

\noindent\textbf{Definition.} A \textit{long virtual knot diagram} is a smooth immersion $f\colon \mathbb{R}\to \mathbb{R}^2$ with crossing information at each double point, and such that $f(x)=(x,0)$ for
$|x|$ sufficiently large.\\
Long knot diagrams are assumed to be oriented from the left to the right.

\noindent\textbf{Definition.} A \textit{long virtual knot} is an equivalence class of long virtual
knot diagrams modulo generalized Reidemeister moves (see \cite{Man} for details).

Long virtual knots and their invariants first appeared in \cite{GPV}.
For the next several definitions we are going to follow the terminology introduced in \cite{FJK}.

\noindent\textbf{Definition.} Let $X$ be a set, and let $P_n(X)$ denote the group of permutations of the $n$-fold Cartesian product, $X^n$.
A \textit{switch} on $X$ is defined as an element $S\in P_2(X)$ satisfying the following relation in $P_3(X)$,
\begin{equation}
(S\times id)(id\times S)(S\times id)=(id\times S)(S\times id)(id\times S). \label{YBE}
\end{equation}
 
Any switch on $X$ gives a representation of the braid group, $B_n$, into the group $P_n(X)$ by sending the standard 
generator of $B_n$, $\sigma_i$, to $S_i=(id)^{i-1}\times S\times (id)^{n-i-1}.$

\noindent\textbf{Example.}
Let $(a,b)\to a*b$ be a rack action on a set $X$. Then $S(a,b)=(b,a*b)$ is a switch, called rack switch. For a thorough introduction to the theory of racks, see \cite{FR}.

\noindent\textbf{Example.}
Let $X=G$ be a group. Then $S(g,h)=(gh^{-1}g^{-1}, gh^2)$ is called the Wada switch (introduced by M. Wada,
see \cite{Wada} and \cite{CESSW}).
As stated in the above definition, any switch $S$ is a bijection.
The inverse to Wada switch is given by $S^{-1}(g,h)=(g^2 h, h^{-1}g^{-1}h).$

\noindent\textbf{Example.}
Let $X$ be a module over a commutative ring. The linear isomorphism $S\colon X^2\to X^2$ given by
$$S(a,b)=(\mu b, \lambda a+ (1-\mu\lambda)b),$$
with $\lambda$, $\mu$ invertible elements of the ring, defines a switch on $X$, called Alexander switch.

A switch $S$ on a set $X$ defines two binary operations on $X$ in the following way 
$$S(a,b)=(b_a,a^b).$$
The operations $b_a$ and $a^b$ are called the down, up operations respectively.
The inverse for $S$ defines two more binary operations called the up-bar and down-bar operations.
They are defined by
$$S^{-1}(a,b)=(b^{\overline{a}},a_{\overline{b}}).$$

The above notation makes brackets obsolete. For example, ${a^b}_c$ means $(a^b)_c$, and $a^{b_c}$ equals to
$a^{(b_c)}$.

The Equation \ref{YBE} leads to the following identities:
\begin{itemize}
\item Up Interchanges: 
\begin{equation}
a^{bc}=a^{{c_b}{b^c}}\label{Upinter}
\end{equation}
\item Down Interchanges:
\begin{equation}
a_{bc}=a_{{c^b}{b_c}}\label{Downinter}
\end{equation}
\item The Rule of Five:
\begin{equation}
{a_b}^{c_{b^a}}={a^c}_{b^{c_a}}\label{Ruleoffive}
\end{equation}
\end{itemize}

We also have the identities called partial inverses:
\begin{equation}
b=b^{\overline{a}a_{\overline{b}}}=b^{a\overline{a_b}}=b_{a\overline{a^b}}=b_{\overline{a}a^{\overline{b}}} \label{pinv}
\end{equation}
for all $a$, $b\in X$.

The up, down, up-bar, and down-bar operations define four endomorphisms of $X$:
\begin{itemize}
\item $x\mapsto x^a$ (also denoted as $f^a$)
\item $x\mapsto x_a$ (also denoted as $f_a$)
\item $x\mapsto x^{\overline{a}}$
\item $x\mapsto x_{\overline{a}}.$
\end{itemize}
\noindent\textbf{Definition.} 
Let $S$ be a switch on $X$. We say that the pair $(X, S)$ defines a \textit{birack} if the following conditions are satisfied:
\begin{enumerate}
\item The map $f^a\colon X\to X$ is a permutation in $P(X)$ for every $a\in X$, and the inverse permutation is denoted as
$x^{a^{-1}}=(f^a)^{-1}(x)$.
\item The map $f_a\colon X\to X$ is a permutation in $P(X)$ for every $a\in X$, and we use the notation 
$x_{a^{-1}}=(f_a)^{-1}(x)$ for the inverse permutation.
\end{enumerate}

Some authors call such structures strong biracks and only require $f^a$ and $f_a$ to be surjective in the definition of birack.
However, most useful biracks that appeared so far in the literature are strong biracks.

The following two lemmas from \cite{FJK} make calculations involving biracks easier.

\begin{lemma}
For a birack, the functions $x\mapsto x^{\overline{a}}$ and $x\mapsto x_{\overline{a}}$ are bijective. The inverses are written
$x\mapsto x^{\overline{a}^{-1}}$ and $x\mapsto x_{\overline{a}^{-1}}$. They are given by the formulae
$$ x^{\overline{a}^{-1}}=x^{a_{x^{-1}}} \ \textrm{and} \ \, x_{\overline{a}^{-1}}=x_{a^{x^{-1}}}.$$
\end{lemma}

\begin{lemma}\label{ids}
Let $x$, $b$, $c$ be elements of a birack $X$. Then the following equalities hold:
\begin{itemize}
\item[(i)] $x_{(b_c)^{-1}}=x_{c^{-1}b^{-1}c^b}$
\item[(ii)] $x_{(c^b)^{-1}}=x_{b_c c^{-1}b^{-1}}$
\item[(iii)] $x^{(b^c)^{-1}}=x^{c^{-1}b^{-1}c_b}$
\item[(iv)] $x^{(c_b)^{-1}}=x^{b^c c^{-1}b^{-1}}.$
\end{itemize}
\end{lemma}

\noindent\textbf{Definition.} We say that the birack $(X, S)$ is a \textit{biquandle} if the following identities hold:
$$ a^{a^{-1}}=a_{a^{a^{-1}}} \ \textrm{and} \ \, a_{a^{-1}}=a^{a_{a^{-1}}}$$
for every $a\in X$.

The Wada switch and Alexander switch are examples of switches that define biquandles.

\section{Biquandle longitude invariant}

In this section we explain how to color virtual link diagrams using biquandles, and we describe our biquandle invariant of long virtual knots.
The details on labeling link diagrams with birack or biquandle elements can be found in \cite{FJK}.
Here, we just present a short description of biquandle colorings.

\begin{figure}
\begin{center}
\includegraphics[height=5.5cm]{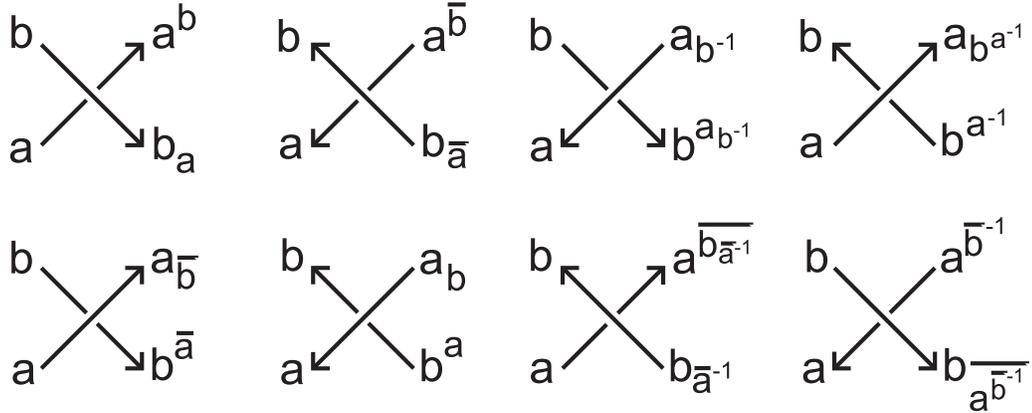}
\caption{Labeling of crossings in biquandle coloring.\label{bicolorings}}
\end{center}
\end{figure}

\begin{figure}
\begin{center}
\includegraphics[height=2.5cm]{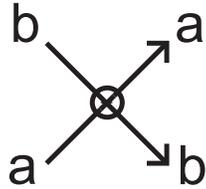}
\caption{Labeling of virtual crossing with biquandle elements.\label{virtcr}}
\end{center}
\end{figure}

\noindent\textbf{Definition.} 
A \textit{semi-arc} of the virtual link diagram is an arc running from a classical crossing to the next classical crossing, ignoring virtual crossings.

\noindent\textbf{Definition.} 
A \textit{biquandle coloring} of the virtual link diagram with a given biquandle $\mathcal{B}$ is an attachment of an element of $\mathcal{B}$ to each semi-arc of the diagram such that at every classical crossing the condition illustrated in 
Figure \ref{bicolorings} is satisfied. 
Notice, that the Figure \ref{bicolorings} shows all possible situations at a crossing with given two labels $a$ and 
$b\in \mathcal{B}$. For a virtual crossing the labels carry across the strings like in Figure \ref{virtcr}.
Alternatively, we can view a biquandle coloring as a function from the set of semi-arcs of the diagram to the biquandle
$\mathcal{B}$, satisfying above conditions.

\begin{theorem}[\cite{FJK}]
Let diagrams $D_1$ and $D_2$ of virtual links be equivalent by a series of generalized Reidemeister moves.
Then any labeling of $D_1$ by elements of a biquandle $\mathcal{B}$ defines a unique labeling of $D_2$ using these moves.
\end{theorem}

It follows that the number of biquandle colorings of a virtual link diagram with a given biquandle $\mathcal{B}$ is an invariant of a virtual link.

To color a diagram of a long virtual knot with elements of some biquandle $\mathcal{B}$, we follow the above procedure for colorings of closed knots. 

Now we proceed to define our invariant.

\noindent\textbf{Definition.} 
Let $D$ be a diagram of a long virtual knot $K$, and let $\mathcal{C}(D)$ denote a biquandle coloring of this diagram with elements of some finite biquandle $\mathcal{B}$.
To every such $\mathcal{C}(D)$ we can assign a bijective map 
$\mathcal{L}(\mathcal{C})\colon \mathcal{B} \to \mathcal{B}$ defined according to the following rules:
\begin{itemize}
\item Let the symbol $x^{u^*}$ denote either $x^u$ or $x^{u^{-1}}$. \\The map $\mathcal{L}(\mathcal{C})$ is of the form $\mathcal{L}(\mathcal{C})(x)=x^{u_1^*u_2^*\ldots u_s^*},$
where $u_1, u_2, \ldots , u_s$ are certain elements of $\mathcal{B}$ encountered in the colored diagram $D$ when traveling along it from the left to the right, and the type of the up-operators $u^*$ is as explained below.
\item Every classical crossing will contribute two up-operators to the definition of $\mathcal{L}$, one when traveling along over-arc of this crossing, and one when traveling along under-arc. Virtual crossings will be ignored.
\item  At every classical crossing of the colored diagram $D$ there are four labels. The label that is going to be used in the definition of $\mathcal{L}$ is the one that is pointed by the normal vector to the traveled path. For a moment, let us name this element as $u$. Normal vectors are chosen so that the pair $($tangent vector, normal vector$)$ matches the usual orientation of the plane.
\item When traveling along under-arc of the positive crossing, we add $u$ to the list of up-operators. If the crossing
is negative, we use $\, u^{-1}$ (notice that it denotes the inverse of operator, not the inverse of an element).
\item When traveling along over-arc, above convention is reversed, i.e.,  we take $u^{-1}$ as an up-operator if the crossing is positive, and $u$ if it is negative.  
\end{itemize}

We call the map $\mathcal{L}(\mathcal{C})$ a \textit{colored biquandle longitude}.

\begin{figure}
\begin{center}
\includegraphics[height=7.5cm]{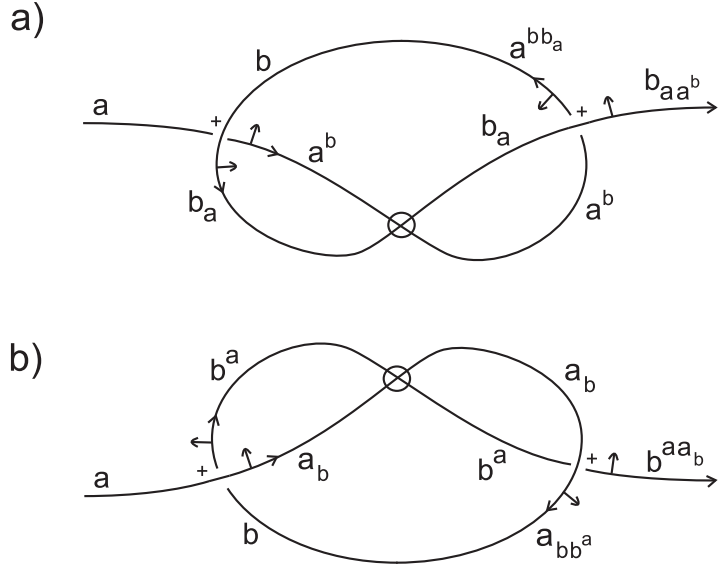}
\caption{Biquandle colorings of the long virtual trefoil knot and its inverse.\label{virttref}}
\end{center}
\end{figure}

\noindent\textbf{Example.} The Figure \ref{virttref}(a) shows an abstract biquandle coloring $\mathcal{C}$ of a long virtual trefoil knot.
Notice that for this coloring to be valid, the condition $b=a^{b b_a}$ has to be satisfied.
For this knot, the map $\mathcal{L}(\mathcal{C})$ is of the form
$ \mathcal{L}(\mathcal{C})(x)=x^{b\, b_a {(a^b)}^{-1}b^{-1}}.$

\begin{theorem} \label{invariant}
The above map $\mathcal{L}(\mathcal{C})$ is invariant under generalized Reidemeister moves.
\end{theorem}

\begin{figure}
\begin{center}
\includegraphics[height=5.2 cm]{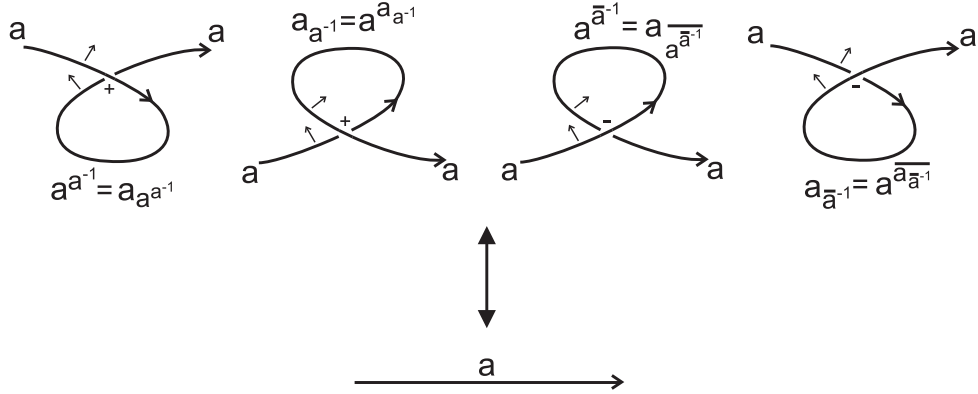}
\caption{Biquandle colorings and the first Reidemeister move.\label{firstrad}}
\end{center}
\end{figure}

\begin{figure}
\begin{center}
\includegraphics[height=5.5 cm]{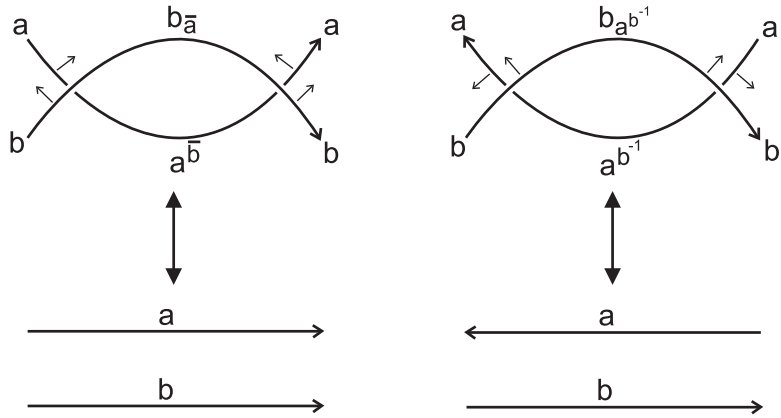}
\caption{Biquandle colorings and the second Reidemeister move.\label{radtwo}}
\end{center}
\end{figure}

\begin{figure}
\begin{center}
\includegraphics[height=5.5cm]{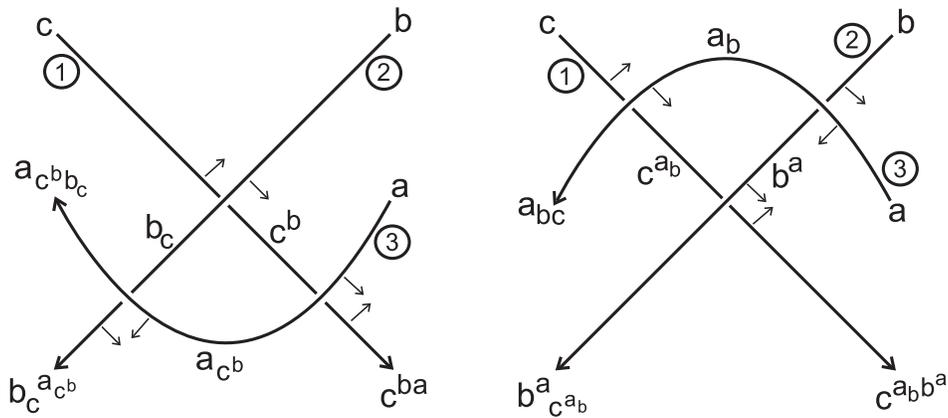}
\caption{Biquandle colorings and the third Reidemeister move.\label{radthree}}
\end{center}
\end{figure}

\begin{proof}
First, let us notice that virtual Reidemeister moves do not affect biquandle colorings. Therefore, they do not change $\mathcal{L}(\mathcal{C})$.

The effect of the first Reidemeister move on biquandle colorings is illustrated in Figure \ref{firstrad}.
In each of these cases the normal vectors to the traveled path point towards two elements that are the same (here, we use conditions from the definition of biquandle). These two elements will be added to the sequence of operators in $\mathcal{L}(\mathcal{C})$ with opposite signs because one is taken when traveling along under-arc, and the other is used when traveling along over-arc.
Thus, they will not influence $\mathcal{L}(\mathcal{C})$.
 
In the case of the second Reidemeister move (Figure \ref{radtwo}), there is a part of colored biquandle longitude that corresponds to the string that is on the top, and the part that is related to the bottom string. Again, when we move along these strings, normal vectors to the traveled path point towards identical elements. Thanks to the fact that the crossing signs are opposite, identical elements will be added to 
$\mathcal{L}(\mathcal{C})$ as opposite operators, and both contributions to the colored biquandle longitude will be trivial.

The proof is not that obvious in the case of the third Reidemeister move (see Figure \ref{radthree}).
We need to consider three parts of $\mathcal{L}(\mathcal{C})$ corresponding
to three strings involved in the third Reidemeister move, and show that they are the same on the left and right sides of the Figure \ref{radthree}. Let us write them as sequences of elements, remembering that they are really sequences of up-operators.
On the left side we have the following parts of $\mathcal{L}(\mathcal{C})$:
\begin{itemize}
\item[(L1)] $b\, a$
\item[(L2)] ${(c^b)}^{-1} a_{c^b}=b^{-1} c^{-1} b_c\, a_{c^b}\ $\\ Above equality holds because of Lemma \ref{ids}(iii).
\item[(L3)] ${(c^{ba})}^{-1}{({b_c}^{a_{c^b}})}^{-1}=a^{-1} {(c^b)}^{-1} a_{c^b}\, {(a_{c^b})}^{-1} {(b_c)}^{-1}
a_{c^b b_c}=$\\ $a^{-1} b^{-1} c^{-1} b_c\, {(b_c)}^{-1} a_{bc}=a^{-1} b^{-1} c^{-1} a_{bc}\ $\\ Here,
we used (three times) Lemma \ref{ids}(iii) and the rule of Down Interchanges (\ref{Downinter}).
\end{itemize}

On the right side of Figure \ref{radthree} the sequences of operators are as follows.
\begin{itemize}
\item[(R1)] $a_b\, b^a= b\,a\ $ (by the rule of Up Interchanges (\ref{Upinter}))
\item[(R2)] $a\, {(c^{ba})}^{-1}=a\, a^{-1} {(c^b)}^{-1} a_{c^b}=b^{-1} c^{-1} b_c\, a_{c^b}$ (by Lemma \ref{ids}(iii))
\item[(R3)] ${(b^a)}^{-1}{(c^{a_b})}^{-1}=a^{-1}b^{-1}a_b\, {(a_b)}^{-1} c^{-1} a_{bc}=
a^{-1}b^{-1}c^{-1}a_{bc}$ (by Lemma \ref{ids}(iii))

\end{itemize}
Corresponding sequences operate in the same way on any biquandle element. It is known that all other types of the third Reidemeister move can be obtained from the move that we considered, and the Reidemeister moves of type one and two. Therefore, $\mathcal{L}(\mathcal{C})$ is invariant under all Reidemeister moves.
\end{proof}

\noindent\textbf{Definition.} Let $D$ be a diagram of a long virtual knot $K$, and let $Col(D,\mathcal{B},p)$ denote the set of colorings of $D$ with biquandle $\mathcal{B},$ satisfying the condition that the label
of the initial arc of $D$ is equal to $p$. Notice that the generalized Reidemeister moves do not change the color of the initial arc.
Therefore, the number of such colorings is a knot invariant. 
We call the family of the corresponding colored biquandle longitudes, 
$$\{\mathcal{L}(\mathcal{C}) \colon \mathcal{C}\in Col(D,\mathcal{B},p)\}$$ a \textit{biquandle longitude invariant} of the knot $K$.

The fact that a biquandle longitude invariant is indeed invariant under generalized Reidemeister moves follows from the
Theorem \ref{invariant}.

One way of utilizing such family of maps is to fix an element $x\in \mathcal{B}$, and consider the formal sum
\begin{equation}
S(D,\mathcal{B},p,x)=\sum_{\mathcal{C}} \mathcal{L}(\mathcal{C})(x), \label{biqsuma}
\end{equation}
taken over all $\mathcal{C}\in Col(D,\mathcal{B},p)$. Naturally, this sum is also a knot invariant. 

\noindent\textbf{Example.} Let us use the above invariant to show that the long virtual trefoil knot from Figure \ref{virttref}(a) is not equivalent to
its inverse (Figure \ref{virttref}(b)).
We are going to use a biquandle defined by Wada switch on a symmetric group $S_5$. 
The eight binary biquandle operations defined by this switch are as follows.
\begin{center}
\begin{tabular}{l  @{\hspace{2.5 cm}} l}
$h_g=gh^{-1}g^{-1}$  &  $h_{g^{-1}}=g^{-1}h^{-1}g$ \\
$g^h=gh^2$  &  $g^{h^{-1}}=gh^{-2}$ \\
$h^{\overline{g}}=g^2h$  &  $h^{{\overline{g}}^{-1}}=g^{-2}h$\\
$g_{\overline{h}}=h^{-1}g^{-1}h$  &  $g_{{\overline{h}}^{-1}}=hg^{-1}h^{-1}$
\end{tabular}
\end{center}
The biquandle longitude for the first knot is
$$ \mathcal{L}_1(\mathcal{C})(x)=x^{b\, b_a {(a^b)}^{-1}b^{-1}}.$$ 
Using the above binary operations we can translate this expression 
into 
$$ \mathcal{L}_1(\mathcal{C})(x)=x b^2 (ab^{-1}a^{-1})^2 (ab^2)^{-2} b^{-2}=x b^2ab^{-2}(a^{-1}b^{-2})^3.$$
The natural choice for $x$, when we use group biquandles, is the group identity. 
As noted before, the condition $b=a^{b b_a}$ is required for the biquandle coloring of the diagram in the Figure \ref{virttref}(a).
With Wada biquandle, this condition can be written as
$$ b=a b^2 (a b^{-1} a^{-1} )^2=a b^2 a b^{-2} a^{-1}.$$
For the inverse of the long virtual trefoil knot (Figure \ref{virttref}(b)), we write the biquandle longitude as
$$ \mathcal{L}_2(\mathcal{C})(x)=x^{{(b^a)}^{-1} {(b^{a a_b})}^{-1} a a_b}.$$
It is equivalent to
$$ \mathcal{L}_2(\mathcal{C})(x)=x a^{-2} b^{-3} a^{-2} b^{-1} a^2 b a^{-2} b^{-1}.$$
The condition for biquandle coloring of this diagram is
$$b=a_{b b^a},$$ that is,
$$b=b a^2 b a b^{-1} a^{-2} b^{-1}.$$ 
Both knots have 240 colorings with the Wada biquandle $S_5$, and both have the same number of biquandle colorings with the fixed color of the initial arc, for any element of $S_5$.  The values of our invariant corresponding to the colorings with the color of the initial arc equal to $(1,2,3,4)$ are as follows:
$$S(D_1,S_5,(1,2,3,4),())=(1,3,2,5,4)+()+(1,5,3,4,2)+(1,4,5,2,3)+(1,2,4,3,5);$$
$$S(D_2,S_5,(1,2,3,4),())=(1,4,2,3,5)+()+(1,2,5,4,3)+(1,3,4,5,2)+(1,5,3,2,4).$$
Here, $()$ denotes the identity permutation.
These two values are different, and therefore, the long virtual trefoil knot is noninvertible.

\end{document}